**UDC 511.3**

**The Generalized Riemann's hypothesis**

S.V. Matnyak

Khmelnytsky, Ukraine

**Summary:** The article presents the proof of the validity of the generalized Riemann's hypothesis on the basis of adjustment and correction of the proof of the Riemann's hypothesis in the work [1], obtained by a finite exponential functional series and finite exponential functional progression.

**Introduction.** In this paper we give the proof of the generalized Riemann's hypothesison the basis of adjustments and corrections to the proof of the Riemann's hypothesis of the zeta function, which was undertaken in [1], as well as values of specified limits of the coefficient $C$. The paper also provides a refutation of the hypothesis of Mertens.

**The formulation of the problem (Generalized Riemann's hypothesis)**. All non-trivial zeros of Dirichlets function $L(s, \chi)$ have a real part that is equal to $\sigma = \dfrac{1}{2}$.

**1).For this we first is will prove the Riemann's hypothesisfor the zeta –function $\varsigma(n)$.**

**The solution**. For the confirmation of the Riemann's hypothesis we will give the definitions and prove the following theorem.

**Definition1**.The expression

$$\sum_{k=1}^{[\sqrt{n}]} u^{\frac{1}{k}}(x) = u(x) + u^{\frac{1}{2}}(x) + u^{\frac{1}{3}}(x) + \ldots + u^{\frac{1}{[\sqrt{n}]}}(x) \tag{1}$$

called finite exponential functional series with respect to the variable exponent $\frac{1}{k}$, where $k = \left\{1, 2, 3, \ldots, \left[\sqrt{n}\right]\right\}$.

**Definition2**.The progression of the type

$$a(x), a(x)q(x), a(x)q^2(x), \ldots, a(x)q^{\sqrt{x}}(x) \tag{2}$$

is called finite  exponential **functional**  progression, if their first member $a(n)$ is  a  function of x or is equal to 1 and the denominator $q(n)$ is a  function of the variable $x^{\frac{1}{\sqrt{x}}}$ .

**Theorem 1**. If the  set of natural numbers  $N_n^+ = \{1,2,...,k,...,n\}$ is the union of subsets $M_1, M_2, M_3, M_4, M_5, M_6$ and  these subsets  are disjoint  and have  appropriately  for  $m, r, s, t, v, u$  the elements, the number of elements of the set  $N_n^+ = M_1 \cup M_2 \cup M_3 \cup M_4 \cup M_5 \cup M$ equals to

$$n = m + r + s + t + v + u$$

.

**Proof**. The theorem is  proved similarly to the theorem 7.11, [4,p.50]

**Theorem 2.** $M(n) = \sum_{k=1}^{\infty} \mu(k)$ of  a series $\frac{1}{\varsigma(s)} = \sum_{n=1}^{\infty} \frac{\mu(n)}{n^s}$ equals. $|M(n)| < 1,5\sqrt{n}$

**Proof.** Let the number  $N = N(n)$ be the quantity of elements of the set of natural numbers $N_n^+ = \{1,2,3,...,n\}$. The positive integers  $N_n^+$ consist: of 1; primes, the quantity of which is denoted by  $\Pi = \Pi(n)$ ; of natural numbers, which are divided on  $p^m$ where quotient  different  from 1 with $m \geq 2$, the number of which is denoted by  $K^\kappa = K^\kappa(n)$; of naturel numbers which are decomposed on product of dual accumulation factors, let is denote the number of these <u>congruous numbers</u> through  $T^n = T^n(n)$, and the amount of numbers, that can be converted into product of unpaired number of primes, will be denoted by  $T'' = T''(n)$. Then the amount of natural numbers  $N$ , of the set of natural numbers  $N_n^+$, according to the theorem 1, equals

$$.N_n^+ = 1 + \Pi + T^n + T'' + K^\kappa + K \qquad (3)$$

The number of the natural series  $N_n^+$ approximately can be expressed as a finite exponential function  series

$$f_1(n) = \sum_{k=2}^{[n]} n^{\frac{1}{k}} = \sqrt{n} + \sqrt[3]{n} + \sqrt[4]{n} + ... + \sqrt[[n]]{n} \qquad (4)$$

We will write the number of natural numbers  $K(n)$, approximately as finite exponential function series consisting of  the first  $(\sqrt{n} - 1)$ members  of  the series (4), we denote it as  $f_2(n)$, then

$$f_2(n) = \sum_{k=2}^{[\sqrt{n}]} n^{\frac{1}{k}} = \sqrt{n} + \sqrt[3]{n} + \sqrt[4]{n} + ... + \sqrt[[\sqrt{n}]]{n} . \qquad (5)$$

In the sum of the series (4) each term of the series is taken, as the amount of natural numbers. For example,

$\sqrt[k]{n} = \left[\sqrt[k]{n}\right]$; $\left[\sqrt{100}\right] = \{1;2;3;4;5;6;7;8;9;10\}$ ; $\left[\sqrt[3]{100}\right] = \{1;2;3;4\}$ and so on.

A series $\left(\sqrt{n} + \sqrt[3]{n} + \sqrt[4]{n} + ... + \sqrt[[n]]{n}\right)$ is greater than the series of the natural number of positive integers $f_1(n) > N_n^+$.

**Definition 3.** The natural numbers that overlap, are called the finite exponential function series of (5), in which they occur more than once. Let that $N = n$ and $\frac{1}{\sqrt{n}} = \frac{1}{\left[\sqrt{n}\right]}$. The function $f_2(n) > K(n)$, because in the function except the numbers $K(n)$ are included the numbers which overlap.

**Definition 4.** Two infinitely great $f(n)$ and $\varphi(n)$, which are not equal to one another $f(n) \neq \varphi(n)$ are called equivalent if $\lim\limits_{n \to \infty} \frac{f(n)}{\varphi(n)} = 1$ when $\varphi(n) \neq 0$.

We assume that the set of natural numbers $M_1 = (1,2,3,...,n-1,n)$ with algebra $A_1 = \langle N, +, \cdot, -, 0, 1 \rangle$ is vector space $P_1$. In this space we set the standard $\|x\|_n = \max\limits_n |n| = n$, and the set of numbers $M_2 = \left(\sqrt[[n]]{n}, \sqrt[(n-1)]{n}, ..., \sqrt{n}\right)$ with algebra $A_2 = \langle \sqrt[[n]]{n}, +, \cdot, -, 0, 1 \rangle$ will be assumed as vector space $P_2$ with standard

$$\|x\|_{\sqrt{n}} = \max\limits_{\sqrt{n}} \left|\sqrt{n}\right| = \sqrt{n} .$$

Then, the denominator of an exponential function of finite progression, which operates in the space $P_1$ will take as $q = n^{\frac{1}{n}}$ , and the denominator of an exponential function of finite progression, which operates in the space $P_2$ , as $q = n^{\frac{1}{\sqrt{n}}}$ .

The finite exponential functional series (4) is approximable by the sum of the finite exponential functional progression

$$\varphi_1\left(n^{\frac{1}{n}}\right) = 1 + n^{\frac{1}{n}} + n^{\frac{2}{n}} + n^{\frac{3}{n}} + ... + n^{\frac{1}{2}} \tag{6}$$

**Proposition 1.**The finite exponential functional series (4) and the sum of the finite exponential functional progression (6) are equivalent.

**Proof.**To prove the equivalence of the finite exponential functional series $f_1(n)$ with the finite functional progression $\varphi_1\left(n^{\frac{1}{n}}\right)$ the sum of the functional series is written in the form of

$$f_1(n) = \left(\sqrt{n} + \sqrt[3]{n} + \ldots + \sqrt[n]{n}\right) < \left(n^{\frac{1}{\sqrt{n}}}\left(n - \sqrt{n} + 1\right) + \left(\sqrt{n} + \sqrt[3]{n} + \ldots + \sqrt[\sqrt{n}]{n}\right)\right),$$

Then let us write that

$$\lim_{n \to \infty} \frac{f_1(n)}{\varphi_1\left(n^{\frac{1}{n}}\right)} = \lim_{n \to \infty} \frac{\left(\sqrt{n} + \sqrt[3]{n} + \ldots + \sqrt[n]{n}\right)}{\left(1 + n^{\frac{1}{n}} + n^{\frac{2}{n}} + n^{\frac{3}{n}} + \ldots + n^{\frac{1}{2}}\right)} = \lim_{n \to \infty} \frac{\sqrt{n} \cdot \left(1 + n^{\frac{1}{3} - \frac{1}{2}} + n^{\frac{1}{4} - \frac{1}{2}} + \ldots + n^{\frac{1}{n} - \frac{1}{2}}\right)}{\sqrt{n} \cdot \left(n^{\frac{-1}{2}} + n^{\frac{1}{n} - \frac{1}{2}} + \ldots + 1\right)} = 1$$

Therefore, in accordance with the definition 4, the functional series and functional progression will be equivalent.

## Proposition 1 is proved.

**Proposition 2**.The finite exponential functional series(4) and the sum of the finite exponential functional progression (6) $2 \cdot S_1(N)$ are equivalent.

**Proof.**The sum ofthe finite exponential functional progression $\varphi_1\left(n^{\frac{1}{n}}\right)$ can be calculated by the formula

$$2 \cdot S_1(N) = 2 \cdot \frac{\sqrt{n} - 1}{n^{\frac{1}{n}} - 1} > 2 \cdot \frac{\sqrt{n} - 1}{\frac{2}{\sqrt{n}}} = \left(n - \sqrt{n}\right) \qquad (7)$$

and then the limit $\lim_{n \to \infty} \frac{f_1(n)}{\varphi_1\left(n^{\frac{1}{n}}\right)}$ will be equal to:

$$\lim_{n\to\infty}\frac{f_1(n)}{\varphi_1\left(n^{\frac{1}{n}}\right)}=\lim_{n\to\infty}\frac{n^{\frac{1}{\sqrt{n}}}\left(n-\sqrt{n}+1\right)+\left(\sqrt{n}+\sqrt[3]{n}+\ldots+\sqrt[\sqrt{n}]{n}\right)}{n-\sqrt{n}}=$$

$$=\lim_{n\to\infty}\frac{\left(n^{\frac{1}{\sqrt{n}}}n-n^{\frac{1}{\sqrt{n}}}\sqrt{n}+n^{\frac{1}{\sqrt{n}}}1\right)+\left(\sqrt{n}+\sqrt[3]{n}+\ldots+\sqrt[\sqrt{n}]{n}\right)}{n-\sqrt{n}}=\lim_{n\to\infty}\frac{n+\sqrt[3]{n}+\sqrt[4]{n}+\ldots+\sqrt[\sqrt{n}]{n}+1}{n-\sqrt{n}}=1,$$

Where $\lim_{n\to\infty}\sqrt[n]{n}=1$ and $0<\sqrt[n]{n}-1<\dfrac{2}{\sqrt{n}}$ [5,p.67].

Therefore, in accordance with the definition 4, the functional series (4) and the finite sum $2\cdot S_1$ of the functional progression (6) at $\left(\sqrt[n]{n}-1\right)\approx\dfrac{2}{\sqrt{n}}$ , when $n\to\infty$ ,will be equivalent.

**Proposition 2 is proved.**

From the expression (4) and (7) one can see that $\varphi_1\left(n^{\frac{1}{n}}\right)$ is within limit of

$$f_1(n)\approx\varphi_1\left(n^{\frac{1}{n}}\right)\ge 2\cdot S_1(n)>n-\sqrt{n}\ .$$

The sum of the finite exponential functional progression (6) with $q=n^{\frac{1}{n}}$ equals $2\cdot S_1(N)\approx\left(n-\sqrt{n}\right)$.When $(n\to\infty),S_1(N)\to\infty$ . We will compare the function $2\cdot S_1(n)$ with the function $f(n)=N_n^+$.We find

$$k_1=\lim_{n\to\infty}\frac{f(n)}{S_1(n)}=\lim_{n\to\infty}\frac{N}{2\cdot\dfrac{\sqrt{n}-1}{\sqrt[n]{n}-1}}<\lim_{n\to\infty}\frac{N}{2\cdot\dfrac{\sqrt{n}-1}{\dfrac{2}{\sqrt{n}}}}=\lim_{n\to\infty}\frac{N}{n-\sqrt{n}}=1.$$

Therefore, $k_1<1$ .

**Lemma 1.** The number of natural numbers that overlap is less than $1,5\sqrt{n}$ .

**Proof.** To prove this proposition let us denote through $K^n$ is the numbers that occur more than once in the finite exponential functional series (4) when $n \to \infty$, and use the exponential functional series $f_2(n)$.

The series (5) is taken is this form to be because it includes all numbers that overlap. This follows from the expression $2^{\sqrt{n}} > N$, $\sqrt{n} \cdot \ln 2 > \ln N$. Two is taken because it is the smallest prime number that can not be decomposed into prime factors.

The finite exponential functional series (5) will be replace by the sum of finite exponential functional progression

$$\varphi_2\left(n^{\frac{1}{\sqrt{n}}}\right) = \left(1 + \sum_{k=1}^{\frac{\sqrt{n}}{2}} n^{\frac{k}{\sqrt{n}}}\right). \tag{8}$$

**Proposition 3.** The finite exponential functional series (5) and the finite exponential functional progression(8) are equivalent.

**Proof.** The functional series (4) can be written as

$$f_1(n) = \sum_{k=\sqrt{n}}^{n} n^{\frac{1}{k}} + \left(n^{\frac{1}{2}} + n^{\frac{1}{3}} + \ldots + n^{\frac{1}{\sqrt{n}}}\right),$$

and the functional progression (6) is as

$$\varphi_1\left(n^{\frac{1}{n}}\right) = \left(1 + n^{\frac{1}{n}} + n^{\frac{2}{n}} + \ldots + n^{\frac{n-2\sqrt{n}}{2n}}\right) + \left(1 + \sum_{k=1}^{\frac{\sqrt{n}}{2}} n^{\frac{k}{\sqrt{n}}}\right).$$

Discard the first terms of the series and progression, we find that

$$f_2(n) = n^{\frac{1}{2}} + n^{\frac{1}{3}} + \ldots + n^{\frac{1}{\sqrt{n}}} \text{, or } \varphi_2\left(n^{\frac{1}{\sqrt{n}}}\right) = \left(1 + \sum_{k=1}^{\frac{\sqrt{n}}{2}} n^{\frac{k}{\sqrt{n}}}\right).$$

We show that

$$f_3(n) = \sum_{k=\sqrt{n}}^{n} n^{\frac{1}{k}} \text{ and } \varphi_3\left(n^{\frac{1}{n}}\right) = \left(1 + n^{\frac{1}{n}} + n^{\frac{2}{n}} + n^{\frac{3}{n}} + \ldots + n^{\frac{n-2\sqrt{n}}{2n}}\right)$$

are equivalent:

$$k_2 = \lim_{n \to \infty} \frac{f_3(n)}{\varphi_3\left(n^{\frac{1}{n}}\right)} = \lim_{n \to \infty} \frac{\sum_{k=\sqrt{n}}^{n} n^{\frac{1}{k}}}{\left(1 + n^{\frac{1}{n}} + n^{\frac{2}{n}} + \ldots + n^{\frac{n-2\sqrt{n}}{2n}}\right)} = \lim_{n \to \infty} \frac{n^{\frac{1}{\sqrt{n}}}\left(n - \sqrt{n} + 1\right)}{2\left(\frac{n-\sqrt{n}}{2}\right)} = 1.$$

It follows that the finite exponential functional series $f_2(n)$ and the finite exponential functional progression $\varphi_2\left(n^{\frac{1}{\sqrt{n}}}\right)$ are equivalent.

**Proposition 3 is proved.**

To prove the theorem, we introduce the functions series

$$f_4\left(n^{\frac{1}{k_1}}\right) = \sum_{k=2k_1}^{\left[\sqrt{n}\right]} n^{\frac{1}{k}} = n^{\frac{1}{4}} + n^{\frac{1}{6}} + \ldots + n^{\frac{1}{\sqrt{n}}} , \qquad (9)$$

where $k_1 = \left\{2, 3, \ldots, \frac{\left[\sqrt{n}\right]}{2}\right\}$

and functional progression

$$\varphi_4\left(n^{\frac{2}{\sqrt{n}}}\right) = \left(1 + n^{\frac{2}{\sqrt{n}}} + n^{\frac{4}{\sqrt{n}}} + \ldots + n^{\frac{\sqrt{n}-8}{4\sqrt{n}}} + n^{\frac{1}{4}}\right). \qquad (10)$$

If we express the series (5), as a series $\left[\sqrt{n}\right] + \left[\sqrt[3]{n}\right] + \left[\sqrt[4]{n}\right] + \ldots + \left[\sqrt[\sqrt[n]{n}]{n}\right]$, than the series (8) is taken $\left[\sqrt{n}\right] + \left[\sqrt[4]{n}\right] + \left[\sqrt[6]{n}\right] + \ldots + \left[\sqrt[2k]{n}\right]$ is such form so that each element of the series (8) overlaps the each element of the series (5) with unpaired exponents of the root. And then we can write that

$$.\left[\sqrt{n}\right] + \left[\sqrt[4]{n}\right] + \left[\sqrt[6]{n}\right] \ldots + \left[\sqrt[k]{n}\right] > \left[\sqrt[3]{n}\right] + \left[\sqrt[5]{n}\right] + \ldots + \left[\sqrt[2k+1]{n}\right]. \qquad (11)$$

Hence the amount of numbers that cover more numbers that overlap.

**Proposition 3 .**The finite exponential functional series (9) and the finite exponential progression (10) are equivalent.

**Proof.** To prove the equivalence of the finite exponential functional series with the finite exponential functional progression in the form of the relation

$$k_3 = \lim_{n \to \infty} \frac{f_4(n)}{\varphi_4(n)} = \lim_{n \to \infty} \frac{n^{\frac{1}{4}} + n^{\frac{1}{6}} + \ldots + n^{\frac{1}{\sqrt{n}}}}{\left(1 + n^{\frac{2}{\sqrt{n}}} + n^{\frac{4}{\sqrt{n}}} + \ldots + n^{\frac{\sqrt{n}-8}{4\sqrt{n}}} + n^{\frac{1}{4}}\right)} = \lim_{n \to \infty} \frac{n^{\frac{1}{4}}\left(1 + n^{\frac{1}{6}-\frac{1}{4}} + \ldots + n^{\frac{1}{\sqrt{n}}-\frac{1}{4}}\right)}{n^{\frac{1}{4}} \cdot \left(n^{-\frac{1}{4}} + n^{\frac{2}{\sqrt{n}}-\frac{1}{4}} + \ldots + 1\right)} = 1$$

Therefore, in accordance with the definition 4, the functional series (9) and a functional progression (10) are equivalent.

**Proposition 3 is proved**.

**Proposition 4.** The finite exponential functional series (9) and the finite sum $4 \cdot S_4(n)$ of the exponential function progression (10) are equivalent when $\sqrt[n]{n} - 1 \approx \frac{2}{\sqrt{n}}$.

**Proof.** The sum of functional series $f_4(n)$ is more than $\left(\frac{\sqrt{n}}{2} - 1\right)$, and the sum of functional progression $\varphi_4\left(n^{\frac{2}{\sqrt{n}}}\right)$ is considered, as the sum of the functional progression with $q = n^{\frac{2}{\sqrt{n}}}$. Then we find that $4 \cdot S_4(n) = 4 \cdot \frac{\sqrt[4]{n}-1}{\sqrt[n]{n^2}-1}$. In order to calculate the function al $\left(n^{\frac{2}{\sqrt{n}}} - 1\right)$, let us set $t^2 = n$ or $t = \sqrt{n}$, and then we obtain

$$\left(n^{\frac{2}{\sqrt{n}}} - 1\right) = \left(t^{\frac{4}{t}} - 1\right) = \left(t^{\frac{2}{t}} - 1\right) \cdot \left(t^{\frac{2}{t}} + 1\right) = \left(t^{\frac{1}{t}} - 1\right) \cdot \left(t^{\frac{1}{t}} + 1\right) \cdot \left(t^{\frac{2}{t}} + 1\right).$$

Since $\lim_{t \to \infty} \sqrt[t]{t} = 1$ then $\left(n^{\frac{2}{\sqrt{n}}} - 1\right) < \frac{8}{\sqrt[4]{n}}$. And then we will have

$$4 \cdot S_4(n) \approx 4 \frac{\sqrt[4]{n}-1}{\frac{8}{\sqrt[4]{n}}} = \frac{\sqrt{n} - \sqrt[4]{n}}{2}, \qquad (12)$$

Using the definition 4 we will have

$$k_3 = \lim_{n \to \infty} \frac{f_4(n)}{\varphi_4\left(n^{\frac{2}{\sqrt{n}}}\right)} = \lim_{n \to \infty} \frac{\frac{\sqrt{n}}{2} - 1}{\frac{\sqrt{n} - \sqrt[4]{n}}{2}} = 1.$$

Therefore, a function of series (9) and the sum $4 \cdot S_4(n)$ of functional progression (10) are equivalent.

**Proposition 4 is proved.**

From the expressions (10) and (12) it is clear that $\varphi_4(n)$ is within

$$f_4\left(n^{\frac{2}{\sqrt{n}}}\right) \approx \varphi_4(n) \ge 4 \cdot S_4(n) > \frac{\sqrt{n} - \sqrt[4]{n}}{2}.$$

Then we compare function $\sqrt{n}$ with the function $f_4(n)$ when $n \to \infty$ and we obtain

$$k_4 = \lim_{n \to \infty} \frac{\sqrt{n}}{4 \cdot \frac{\sqrt[4]{n} - 1}{\sqrt[n]{n^2} - 1}} < \lim_{n \to \infty} \frac{\sqrt{n}}{4 \cdot \frac{\sqrt[4]{n} - 1}{\frac{8}{\sqrt[4]{n}}}} = \lim_{n \to \infty} \frac{8 \cdot \sqrt{n}}{4 \cdot \left(\sqrt{n} - \sqrt[4]{n}\right)} = 2.$$

Hence, we have that $k_4 < 2$, or $\sqrt{n} - 2 \cdot f_4\left(n^{\frac{2}{\sqrt{n}}}\right) < 0$.

Therefore, $\frac{\sqrt{n}}{2} < f_4(n)$. Let us take into account the value of the finite exponential functional series $f_4(n)$, and write that $\frac{\sqrt{n}}{2} < K^n(n) < \sqrt{n} + \frac{\sqrt{n}}{2} = \frac{3}{2}\sqrt{n}$.

**Lemma 1 is proved.**

Then we can write that

$$N(n) - C \cdot f_1(n) < 0, \qquad \text{(13)}$$

Using the inequality $f_1(n) > 2 \cdot K(n)$ we will write that

$$f_1(n) \approx K(n) + \left(n - \sqrt{n} + 1\right) + K^n(n) \qquad \text{(14)}$$

If we substitute value of the function $f_2(n)$ (14) into (13), we obtain

$$N(n) - C \cdot \left(K(n) + \left(n - \sqrt{n} + 1\right) + K^n(n)\right) < 0.$$

Using Lemma 1, we obtain

$$N(n) - C \cdot \left( K(n) + \left( n - \sqrt{n} + 1 \right) + \frac{3}{2} \sqrt{n} \right) < 0.$$

Hence; we find that

$$N(n) - C \cdot \left( K(n) + \left( n - \sqrt{n} + 1 \right) \right) < C \cdot \frac{3}{2} \sqrt{n}.$$

The value $N$ from (3) is substituted instead $N$, we obtain

$$1 + \Pi + T^H + T^n + K^K + K(n) < C \cdot \left( 1{,}5\sqrt{n} + K(n) + \left( n - \sqrt{n} + 1 \right) \right), \text{ or}$$

$$1 + \Pi + T^H + T^n + K^K < C \cdot \left( 1{,}5\sqrt{n} + \left( 1 + \Pi + T^H + T^n + K^K - \sqrt{n} + 1 \right) \right). \qquad (15)$$

Then we can write that appropriately of the properties of the function of Mobius- $\mu(n) = 1$ ,when $n = 1$; $\mu(n) = (-1)^k$ ,where $k$ is the amount of prime factors of the numbers $n = p_1 \cdot p_2 \cdot ... \cdot p_k$ and $\mu(n) = 0$ when $n$ is multiple $p^m$ for $m \geq 2$,

$$1 + T^n - \left( \Pi + T^H \right) < C \cdot \left( \left( 1 + T^n \right) - \left( \Pi + T^H \right) - \sqrt{n} + 1 + 1{,}5\sqrt{n} \right). \qquad (16)$$

We write that

$$M(n) = 1 + T^n - \left( \Pi + T^H \right).$$

Then the expression(16) takes the form

$$M(n) < C \cdot \left( \left( M(n) - \sqrt{n} + 1 \right) + 1{,}5\sqrt{n} \right).$$

Therefore

$$\left| M(n) \right| < \frac{C}{C - 1} \cdot \left( 0{,}5\sqrt{n} + 1 \right). \qquad (17)$$

**The theorem 2 is proved.**

**2) For the Mertens function you can find a more  precise estimate.**

   **Lemma 2**.The accurate assessment.   $M(n) = \sum_{k=1}^{\infty} \mu(k)$ in a series $\dfrac{1}{\varsigma(s)} = \sum_{n=1}^{\infty} \dfrac{\mu(n)}{n^s}$ will be equal

$$M(n) < \frac{C}{C-1}\left(1{,}25 \cdot \sqrt{n} + \sqrt[3]{n} - 1{,}25 \cdot \sqrt[4]{n} - \sqrt{n} + 1\right).$$

**Proof**. In order to find a more accurate estimate than $M(n) < \frac{C}{C-1}\left(0{,}5\sqrt{n}+1\right)$, let us find the sum of the finite exponential functional series (5) $f_2(n) = \sum\limits_{k=2}^{\left[\sqrt[4]{n}\right]} n^{\frac{1}{k}} = \left(\sqrt{n} + \sqrt[3]{n} + \sqrt[4]{n} + \ldots + \left[\sqrt[4]{n}\right]\sqrt{n}\right)$.

For that we use the functional progression (6) $\varphi_2\left(n^{\frac{1}{\sqrt{n}}}\right) = \left(1 + \sum\limits_{k=1}^{\frac{\sqrt{n}}{2}} n^{\frac{k}{\sqrt{n}}}\right)$, then we obtain that

$$f_2(n) = \sqrt{n} + \sqrt[3]{n} + \sqrt[4]{n} + \ldots + \sqrt[\sqrt{n}]{n} = 2{,}5\sqrt{n} + 2\sqrt[3]{n} - 2{,}5 \cdot \sqrt[4]{n} \quad.$$

We find from the expression (10) that the quantity of numbers that overlap is less than

$\frac{f_2(n)}{2}$ because $\frac{f_2(n)}{2} > \left[\sqrt[3]{n}\right] + \left[\sqrt[5]{n}\right] + \ldots + \left[\sqrt[2k+1]{n}\right]$. Using the method given determine

$M(n) < \frac{C}{C-1}\left(0{,}5 \cdot \sqrt{n} + 1\right)$ and the expression $\frac{f_2(n)}{2} > \left[\sqrt[3]{n}\right] + \left[\sqrt[5]{n}\right] + \ldots + \left[\sqrt[2k+1]{n}\right]$ we obtain that

$$M(n) < \frac{C}{C-1}\left(1{,}25 \cdot \sqrt{n} + \sqrt[3]{n} - 1{,}25 \cdot \sqrt[4]{n} - \sqrt{n} + 1\right) \quad.$$

Hence, we find that the upper limit of the value functions $\lim\limits_{n\to\infty} \frac{M(n)}{\sqrt{n}}$ will be the value

$$\lim\limits_{n\to\infty} \sup \frac{M(n)}{\sqrt{n}} < \frac{C}{C-1} \cdot 0{,}25 \quad,$$

and the lower limit is

$$\lim\limits_{n\to\infty} \inf \frac{M(n)}{\sqrt{n}} > -\frac{C}{C-1} 0{,}25 \quad.$$

Therefore, the evaluation $M(n) < \frac{C}{C-1}\left(0{,}25\sqrt{n} + \sqrt[3]{n} - 1{,}25\sqrt[4]{n} + 1\right)$ is a more accurate estimate than $M(n) < \frac{C}{C-1}\left(0{,}5 \cdot \sqrt{n} + 1\right)$ when $n \to \infty$.

**Lemma 2 is proved.**

3) The theorem 2 proves that the upper limit value of the function $\lim\limits_{n\to\infty}\dfrac{M(n)}{\sqrt{n}}$ equals

$$\lim_{n\to\infty}\sup\frac{M(n)}{\sqrt{n}} < 0,5\cdot\frac{C}{C-1},$$

and the lower limit is

$$\lim_{n\to\infty}\inf\frac{M(n)}{\sqrt{n}} > -0,5\cdot\frac{C}{C-1},$$

**Proposition 4.** $\dfrac{C}{C-1}\cdot\left(0,5\sqrt{n}+1\right) << n^{\frac{1}{2}+\varepsilon}$ when $n\to\infty$

**Proof.** According to the theorem 54 [7.p.114] we have that $M(n)=O\left(n^{\frac{1}{2}+\varepsilon}\right)$. The value

$M(n)=O\left(\dfrac{C}{C-1}\cdot\left(0,5\sqrt{n}+1\right)\right)$ is compared with $M(n)=O\left(n^{\frac{1}{2}+\varepsilon}\right)$, we will write that

$\left(\dfrac{C}{C-1}\left(0,5\sqrt{n}+1\right)\right)=n^{\frac{1}{2}+\varepsilon_0}$ when $n\to\infty$. Hence we find that $\varepsilon_0=\dfrac{\ln\left(\dfrac{C}{C-1}\left(0,5\sqrt{n}+1\right)\right)}{\ln n}$ when

$n\to\infty$, $\varepsilon_0\to 0$. Therefore, we can assume that $\varepsilon>\varepsilon_0$, where $\varepsilon$ is a random small number. And here we find that $\dfrac{C}{C-1}\cdot\left(0,5\sqrt{n}+1\right) << \left(n^{\frac{1}{2}+\varepsilon}\right)$ when $n\to\infty$.

**Proposition 4 is proved.**

**Theorem 3.** The series $\dfrac{1}{\varsigma(s)}=\sum\limits_{n=1}^{\infty}\dfrac{\mu(n)}{n^s}$ converges if $\sigma=\dfrac{1}{2}+\varepsilon>\dfrac{1}{2}$ and

$|M(n)|<\dfrac{C}{C-1}\left(0,5\sqrt{n}+1\right)$ where $\varepsilon$ is a random small number.

**Corollary of Theorem 3 (the Riemann's hypothesis)**. All non-trivial zeros of the

zeta−function have a real part equal to $\sigma=\dfrac{1}{2}$.

**Proof.** A necessary and sufficient condition for the validity of the Riemann's hypothesis is the convergence of the series $\dfrac{1}{\varsigma(s)}=\sum\limits_{n=1}^{\infty}\dfrac{\mu(n)}{n^s}$ when $\sigma>\dfrac{1}{2}$ [7,p.114]. We find the convergence of

the series, when

$$M(n)=O\left(\frac{C}{C-1}\left(0,5\sqrt{n}+1\right)\right) \text{ and } \sigma=\frac{1}{2}.$$

$$\frac{1}{\xi\left(\frac{1}{2}\right)} = \sum_{n=1}^{\infty} \frac{M(n)}{n^{\frac{1}{2}}} = \sum_{n=1}^{\infty} \frac{M(n)-M(n-1)}{n^{\frac{1}{2}}} = \sum_{n=1}^{\infty} M(n) \cdot \left(\frac{1}{n^{\frac{1}{2}}} - \frac{1}{(n+1)^{\frac{1}{2}}}\right) =$$

$$= \sum_{n=1}^{\infty} M(n) \left(\frac{(n+1)^{\frac{1}{2}} - n^{\frac{1}{2}}}{n^{\frac{1}{2}} \cdot (n+1)^{\frac{1}{2}}}\right) \le \sum_{n=1}^{\infty} \frac{M(n)}{2n\sqrt{n}} = \sum_{n=1}^{\infty} \frac{M(n)}{2n\sqrt{n}} = \sum_{n=1}^{\infty} \frac{\frac{C}{C-1}\left(0,5\sqrt{n}+1\right)}{2n\sqrt{n}} = \frac{0,5\frac{C}{C-1}}{2}\sum_{n=1}^{\infty} \frac{1}{n} \to \infty$$

the series diverges.

And when $\sigma = \frac{1}{2} + \varepsilon > \frac{1}{2}$ we have

$$\frac{1}{\varsigma\left(\frac{1}{2}+\varepsilon\right)} = \sum_{n=1}^{\infty} \frac{M(n)}{n^{\frac{1}{2}+\varepsilon}} \le \frac{0,5 \cdot \frac{C}{C-1}}{2} \sum_{n=1}^{\infty} \frac{1}{n^{1+\varepsilon}} \le \frac{0,5 \cdot \frac{C}{C-1}}{2} \int_{1}^{\infty} \frac{1}{n^{1+\varepsilon}} \, dn = \frac{0,5 \cdot \frac{C}{C-1}}{2\varepsilon} = 0,25 \cdot \frac{C}{\varepsilon \cdot (C-1)}$$

the series converges, where ε is an arbitrary small number.

Therefore, the series $\frac{1}{\zeta(s)} = \sum_{n=1}^{\infty} \frac{\mu(n)}{n^s}$ converges uniformly for $\sigma = \frac{1}{2} + \varepsilon > \frac{1}{2}$, and since it is a function $\frac{1}{\varsigma(s)}$ if $\sigma > 1$, for the theorem of analytic continuation, it is also at its $\frac{1}{2} < \sigma \le 1$.

Therefore, the Riemann's hypothesis is true.

**The theorem 3 is proved.**

**4) A determination the values of coefficient $C$.**

Then we can write that according to the properties of Mobius function - $\mu(n)=1$, then $n=1$;

$\mu(n)=(-1)^k$, where $k$ the number of prime factors of the number $n = p_1 \cdot p_2 \cdot ... \cdot p_k$ and $\mu(n)=0$,

when $. n$ is the multiple of $p^m$ for $m \ge 2$ that

$$\left(\left(1+T^n\right)-\left(\varPi+T^{\mu}\right)\right)-C \cdot \left(1+T^n-\left(\varPi+T^{\mu}\right)\right) < \left(0,5 \cdot \sqrt{n}+1\right),$$

Then

$$M(n) < \frac{C}{1-C} \cdot \left(0,5 \cdot \sqrt{n} + 1\right); \left|M(n)\right| < \frac{C}{C-1} \cdot \left(0,5 \cdot \sqrt{n} + 1\right).$$

From the expression $K(n) + \left(n - \sqrt{n} + 1\right) + \frac{3}{2} \cdot \sqrt{n} > 0$, using the properties of Möbius function, it can be written that

$$M(n) < 1,5 \cdot \sqrt{n} - 1.$$

And from the expression $N - \sqrt{n} > 0$ we find that

$$M(n) > \sqrt{n}.$$

This coincides with the results [3]. Then we can find the extent to which the coefficient $C$ is located. From the double inequality $\sqrt{n} < \frac{C}{C-1} \cdot \left(0,5 \cdot \sqrt{n} + 1\right) < 1,5 \cdot \sqrt{n}$, we find that $2 < \frac{C}{C-1} < 3$. And here we find that $1,5 < C < 2$.

.   Using a more precise value $M(n)$, we find that

$$\sqrt{n} < \frac{C}{C-1} \left(0,25\sqrt{n} + \sqrt[3]{n} - 1,25\sqrt[4]{n} + 1\right) < 1,25 \cdot \sqrt{n}.$$

And here we find that $4 < \frac{C}{C-1} < 5$ and from the double inequality we find that the coefficient $c$ will be in the range $1,25 < C < 1,33$

**5) Theorem 4. All non-trivial zeros of Dirichlets function $L(s, \chi)$ have a real part that is equal to $\sigma = \frac{1}{2}$.**

**Proof.** Let's consider the Dirichlet's series

$$L(s, \chi) = \sum_{n=1}^{\infty} \frac{\chi(n)}{n^s}, \, s = \sigma + it, \tag{18}$$

where $\chi$ is the character of modulus $m$.

There is $\varphi(m)$ of such series where $\varphi$ is the Euler's function. Since $\left|\chi(n)\right| \le 1$, the series (18) converges when $\sigma > 1$, as can be seen from a comparison of this series with the series $\sum \dfrac{1}{n^s}$ . We denote it by the sum through series $L(1, \chi)$ . For various characters $\chi$ , we obtain different functions $L(s, \chi)$ .They are called $L$ is the Dirichlet's functions. In studying the properties of these functions it is convenient to distinguish the cases where $\chi$ is the main character $\chi_1$ and when

$\chi \neq \chi_1$ .

a) If $\chi \neq \chi_1$ than the series (18) converges in the semiplane $\sigma > 0$ . Let us show from the

beginning, that the partial sums $\sum\limits_{n<x} \chi(n)$ are limited.We divide the integer number from

1 to $[x]$ into classes of deductions by $\mod m$ and write $[x] = m \cdot q + r$ , $0 \leq r \leq m-1$ . Then

$$\sum_{n<x} \chi(n) = \sum_{n=1}^{[x]} \chi(n) = \left( \sum_{1}^{m} + \sum_{m+1}^{2m} + ... + \sum_{m(n-1)+1}^{mq} \right) \cdot \chi(n) + \sum_{mq+1}^{mq+r} \chi(n).$$

Because of the orthogonality relations

$$\sum_{n(\mod m)} \chi(n) = \begin{cases} \varphi(m), than \chi = \chi_1 \\ 0, than \chi \neq \chi_1 \end{cases}$$

we have

$$\sum_{n<x} \chi(n) = \sum_{mq+1}^{mq+r} \chi(n),$$

hence

$$\left| \sum_{n<x} \chi(n) \right| \leq \sum_{mq+1}^{mq+r} |\chi(n) \leq r < m|.$$

Since $n^{-\sigma}$ at $\sigma > 0$ decreases monotonically and tends to zero when $n \to \infty$ , then the series $\sum \chi(n)/n^s$ converges for real $s = \sigma > 0$ , and, consequently,for all $s$ in the semiplane

$\sigma > 0$ when $\chi \neq \chi$ . If, however, $\sigma < 0$ , then this the series obviously diverge. It's abscissa

converges $\sigma_0 = 0$ and the abscissa of absolute convergence $\bar{\sigma} = 1$ . By the theorem 4, "The

Dirichlet's series $\sum\limits_{n=1}^{\infty} a_n \cdot n^{-s}$ in the semiplane of the convergence is a regular analytic function from $s$ , the successive derivatives of which are obtained by the term differentiation of this the series [8, p.153],the function $L(s, \chi)$, $\chi \neq \chi$ is a regular analytic function from $s$ when $\sigma > 0$.

b) If $\chi = \chi_1$ we use

$$L(s, \chi) = \sum_{n=1}^{\infty} \frac{\mu(n) \cdot \chi(n)}{n^s}, \ \text{Re} \ s \geq \frac{1}{2} . \tag{19}$$

From the theorem 3 it follows that $L(s, \chi) \neq 0$ when $\sigma \geq \frac{1}{2}$. If $\chi_1$ is the main character by $\mod m$, then

$$\chi_1(a) = \begin{cases} 1, then(a, m) = 1 \\ 0, then(a, m) > 1 \end{cases} .$$

Using the condition $\left| \chi(n) \cdot n^{-s} \right| \leq n^{-s}$ the function (19) can be written as

$$L(s, \chi) = \sum_{n=1}^{\infty} \frac{\mu(n) \chi(n)}{n^s} = \sum_{n=1}^{\infty} \frac{\mu(n)}{n^s}, \text{ when } \chi(n) = 1 \text{ and } \sigma \geq \frac{1}{2} .$$

Using the results of the theorem 3, it can be argued that the generalized Riemann's hypothesis is true, and accordingly to it: "All non-trivial zeros of the Dirichlet's functions have a real part equal to $\sigma = \frac{1}{2}$ ".

**The theorem 4 is proved.**

### Appendix

**The proof of the Riemann's hypothesis and the conjecture of Birch and**

**Suinertonna-Daeyr without coefficient $C$.**

A better proof of the Riemann's hypothesis and the proof of the conjecture of Birch and Suinertonna-Daeyr and as a refutation of the Mertens hypothesis, can be found on the basis of this work and when we take the series of

$$f_1 = \left( \sqrt{n} + \sqrt[3]{n} + ... + \sqrt[n]{n} \right) < \left( n^{\frac{1}{\sqrt{n}}} \cdot \left( n - \sqrt{n} + 1 \right) + \left( \sqrt{n} + \sqrt[3]{n} + ... \sqrt[\sqrt{n}]{n} \right) \right)$$

and it can be written in the form

$$\left(n^{\frac{1}{\sqrt{n}}}\cdot\left(n-\sqrt{n}+1\right)+\left(\sqrt{n}+\sqrt[3]{n}+...\sqrt[\sqrt{n}]{n}\right)\right)<n+2,5\cdot\sqrt[2]{n}+2\cdot\sqrt[3]{n}-2,5\cdot\sqrt[4]{n}-\sqrt{n}+1=$$

$$=n+1,5\cdot\sqrt{n}+2\cdot\sqrt[3]{n}-2,5\cdot\sqrt[4]{n}+1 \ ,$$

where

$$f_2\left(n\right)=\sqrt{n}+\sqrt[3]{n}+...+\sqrt[\sqrt{n}]{n}<2,5\sqrt{n}+2\cdot\sqrt[3]{n}-2,5\cdot\sqrt[4]{n} \ .$$

Then we can write

$$N+2,5\cdot\sqrt{n}+2\cdot\sqrt[3]{n}-2,5\cdot\sqrt[4]{n}-\sqrt{n}+1>0$$

or

$$N>-1,5\cdot\sqrt{n}-2\cdot\sqrt[3]{n}+2,5\cdot\sqrt[4]{n}-1 \tag{20}$$

The value $N$ from the expression 3 is substituted instead of $N$ , we obtain

$$N=1+\Pi+T^n+T^u+K^\kappa+K\left(n\right)$$

From the expression (20) we obtain

$$1+\Pi+T^n+T^u+K^\kappa+K\left(n\right)>-1,5\cdot\sqrt{n}-2\cdot\sqrt[3]{n}+2,5\cdot\sqrt[4]{n}-1 \tag{21}$$

The properties of Mobius function [3, p. 3] will be applied to the expression (21) and we obtain that

$$M\left(n\right)>-1,5\cdot\sqrt{n}-2\cdot\sqrt[3]{n}+2,5\cdot\sqrt[4]{n}-1$$

or

$$M\left(n\right)<1,5\cdot\sqrt{n}+2\cdot\sqrt[3]{n}-2,5\cdot\sqrt[4]{n}+1.$$

It will be the smallest value of the function of Mertens $M\left(n\right)$ and the biggest value for the function of Mertens $M\left(n\right)$. From the expression 6

$$\varphi_1\left(n^{\frac{1}{n}}\right) = 1 + n^{\frac{1}{n}} + \ldots + n^{\frac{1}{2}}$$

we find that

$$N - \sqrt{n} > 0 \ , \qquad\qquad (22)$$

when

$$\left(\sqrt[n]{n} - 1\right) \approx \frac{2}{\sqrt{n}} \ \text{ and } \ n \to \infty \ .$$

We write the expression 22 in the form

$$1 + \Pi + T^n + T^u + K^\kappa + K(n) > \sqrt{n} \ . \qquad\qquad (23)$$

Let us apply the properties of Mobius function to the expression 23 and we obtain

$$M(n) > \sqrt{n} \ .$$

Then we can state that the function of Mertens $M(n)$ is within

$$\sqrt{n} < M(n) < 1{,}5 \cdot \sqrt{n} + 2 \cdot \sqrt[3]{n} - 2{,}5 \cdot \sqrt[4]{n} + 1$$

And it rejects the hypothesis of Mertens. And the Riemann's hypothesis respectively by the theorem 3 is the true. And then coefficient $k_1$ in the paper [2, p. 3] will be equal $\approx 2{,}67$. What clarifies the proof of the conjecture of Birch and Suinertona-Dyer.

July 22, 2015

### References

[1] S.V.Matnyak . Proof of the Riemann's hypothesis .// arXiv:1404.5872 [math.GM],22 Apr 2014.

[2] S.V.Matnyak.The proof of the correctness of the Birch and Swinnerton–Diyer conjecture.//arXiv:1406.2270 [math.GM],6 Jun 2014.

[3].A.M. Odlyzko and Herman te Riele. Disproof of the Mertens Conjeture. Jornal fur die reine und angewandte Mathematik. 357. (1985) pp. 138-160.

[4]'. Y.S. Lyapin, A.Y. Yevseyev. Algebra and theory of numbers, p.1. Numbers. Educational book   for students of the faculties of physics  and  mathematics  of teaching colleges. M.: " Prosveshcheniye", 1974,– 383 p.

[5]. G. M. Fihtengolc. Differential and integral calculus course. V. 1.M.:"Nauka", 1969, – 607 p.

[6]. L.Y. Kulikov Algebra and theory of numbers. Work-book for teaching colleges.– M.:Vyshaya shkola, 1979 .– 559 p

[7].Y. K. Titchmarsh. Riemann zetafunction. M.: IL, 1947, –154 p.

[8].K. Chandrasekharan. Introduction toanalytic number theory. M.: Mir,1974, –187p.

E-mail:matn111@yandex.ua